\documentclass[9pt,amstex]{article}
\usepackage{amssymb,amsmath} 
\usepackage{latexsym}

\newtheorem{dfn}{Definition}[section] 
\newtheorem{rmk}{Remark}[section]
\newtheorem{thm}{Theorem}[section] 
\newtheorem{cor}{Corollary}[section]
\newtheorem{prop}{Proposition}[section] 
\newtheorem{lem}{Lemma}[section]

\def\cyclic{\mathop{\kern0.9ex{{+}
\kern-2.2ex\raise-.28ex\hbox{\Large\hbox
{$\circlearrowright$}}}}\limits}

\def\buildrel#1_#2^#3{\mathrel{\mathop{\kern 0pt#1}\limits_{#2}^{#3}}}

\newcommand{\Pf}{{\em Proof}. }
\newcommand{\EPf}
{%
\mbox{}%
\nolinebreak%
\hfill%
\rule{2mm}{2mm}%
\medbreak%
\par%
}

\newcommand{\id}{\mbox{$\mathtt{Id}$}}

\newcommand{\Ad}{\mbox{$\mathtt{Ad}$}}

\newcommand{\C}{\mathbb C}

\newcommand{\B}{\mathbb B} 
\renewcommand{\S}{\mathbb S}

\newcommand{\R}{\mathbb R}

\newcommand{\N}{\mathbb N}

\newcommand{\g}{{\mathfrak{g}}{}}

\renewcommand{\k}{{\mathfrak{k}}{}} 
\newcommand{\p}{{\mathfrak{p}}{}} 
\newcommand{\Der}{{\mathfrak{Der}}{}}

\newcommand{\s}{{\mathfrak{s}}{}} 
\newcommand{\ddto}{{\frac{d}{dt}|_{0}}{}}

\newcommand{\h}{{\mathfrak{h}}{}} 
\renewcommand{\sp}{{\mathfrak{sp}}{}} 
\renewcommand{\a}{{\mathfrak{a}}{}} 
\renewcommand{\b}{{\mathfrak{b}}{}} 
\newcommand{\z}{{\mathfrak{z}}{}}
\newcommand{\CO}{{\cal O}{}}
\newcommand{\CM}{{\cal M}{}}

\newcommand{\CP}{\mathcal P}
\newcommand{\CS}{\mathcal S}
\newcommand{\CE}{\mathcal E}

\newcommand{\CH}{\mathcal H}

\newcommand{\CU}{\mathcal U}
\newcommand{\CF}{\mathcal F}

\def\cref#1{Corollary~\ref{#1}}

\title{Non-formal deformation quantizations of solvable Ricci-type symplectic symmetric spaces}

\author{
{\bf Pierre Bieliavsky}\\
Universit\'e Catholique de Louvain, Belgium.\\
e-mail: {\tt bieliavsky@math.ucl.ac.be}\\
}
\date{} 

\begin{document}

\maketitle
\begin{abstract}
Ricci-type symplectic manifolds have been introduced and extensively studied by M. Cahen et al. \cite{CB,CGS}. In this note, we describe their deformation quantizations in the split solvable  symmetric case. In particular, we
introduce the notion of non-formal {\sl tempered} deformation quantization on such a space. We show that
the set of tempered deformation quantizations is in one-to-one correspondence with the space of Schwartz
operator multipliers on the real line. Moreover we prove that every invariant formal star product on
a split Ricci-type solvable symmetric space is an asymptotic expansion of a tempered non-formal quantization.
This note  illustrates and partially reviews  through an example a problematic studied by the author regarding
non-formal quantization in presence of large groups of symmetries. 
\end{abstract}

\section{Introduction}
A formal universal deformation formula (or Drinfeld twist) based on the enveloping algebra $\CU(\b)$ of a Lie algebra $\b$
is the data of an element $F=\sum_k\theta^kF_k\in\CU(\b)[[\theta]]\otimes\CU(\b)[[\theta]]$ satisfying certain cocycle conditions
(see e.g. \cite{GZ}) equivalent to requiring that the following power series in the formal parameter $\theta$:
\begin{equation*}
\tilde{\star}\;:=\;\sum_k\theta^k\tilde{F_k}
\end {equation*}
is a left-invariant formal star product on the connected simply connected Lie group $\B$ whose Lie algebra is $\b$
(in the above formula we denote by $\tilde{P}$ the left-invariant differential operator on $\B$ associated with the element $P$ of $\CU(\b)$). The symplectic leave $\S$ through the unit element $e$ of $\B$ associated 
with the Poisson structure $[F_1]$ is seen to be a Lie subgroup of $\B$. A  first step in the study of Drinfeld
twists is therefore to consider left-invariant star products on {\sl symplectic Lie groups} (i.e. Lie groups endowed with
left-invariant symplectic structures)\footnote{The abelian case has been first---and beautifully--- treated by Rieffel in the 
$C^\star$ algebraic case \cite{Ri}.}. An important class of the later groups is constituted by Iwasawa subgroups
of Hermitean type simple Lie groups (a simple Lie group $G_s$ is called of {\sl Hermitean type} is the 
center $\z(\k_s)$ of the Lie algebra $\k_s$ of one of its maximal compact subgroups $K_s$ is non trivial). When the homogeneous space $G_s/K_s$ has rank one, the corresponding Iwasawa subgroup $\S$ is a solvable extension
of an Heisenberg group. This is the latter particular situation attended in this note.

\noindent It turns out that the solvable Lie group $\S$ underlies a specific symmetric affine geometry. Of course,
the Iwasawa decomposition $G_s=\S K_s$ tells us that $\S$ may be viewed as a model for the Hermitean
symmetric space $G_s/K_s$. But it is rather a `contaction' of the latter that we're going to consider here (see Section \ref{RICCI}). Any quantization of this `contracted geometry' will yield a left-invariant star product on $\S$, hence a Drinfel'd twist. Before this, we first recall basics on symmetric space theory.

\begin{dfn}A {\sl symmetric space} is a pair $(M,s )$ where $M$ is a
connected smooth manifold and $s:M\times M\rightarrow M$ is a smooth map such
that 
\begin{enumerate}	
\item[(i)] For any $x$ in $M$, $s_x:M \to M:y \to s(x,y)$ is a
involutive diffeomorphism of $M$,	which	admits	$x$ as isolated fixed point. The map $s_x$ is	called	the
{\sl symmetry} at $x$. 
\item[(ii)] For	any $x,y$ in $M$ one has~:
$$s_xs_ys_x=s_{s_xy}\;.$$ 
\end{enumerate} 
Two such spaces
$(M_i,s^{(i)})\,(i=1,2)$ are said  {\bf isomorphic} if there exists a
 diffeomorphism $\varphi :M_1\rightarrow M_2$	such that for any $x_1$ in
$M_1$, 
$$\varphi \circ s^{(1)}_{x_1}=s^{(2)}_{\varphi	 (x_1)}\circ \varphi\;.$$
\end{dfn}

\begin{dfn} 	A	involutive Lie algebra $c=(\g,\sigma )$, ---i.e.
 $\g$ is a finite dimensional real Lie algebra and $\sigma$ is an involutive
automorphism of $\g$ ---is	called	a {\sl symmetric couple} if the
following properties are satisfied: 
\begin{enumerate}	
\item[(i)]	Let $\g=\k\oplus \p$   where $\k$ (resp. $\p$)
is the $+1$ (resp. $-1$) eigenspace of $\sigma$, then 
$$[\p,\p]=\k\;.$$

\item[(ii)] The representation of $\k$ on $\p$, given by the	adjoint
action,	is	faithful. 
\end{enumerate} The dimension	of $\p$ defines	the
{\sl dimension}	of	the triple. Two	such pairs  $(\g_i,\sigma_i)$ $(i=1,2)$ are said	{\sl isomorphic}	if
there exists a Lie algebra isomorphism $\psi :\g_1\rightarrow\g_2$	such
that $\psi\circ\sigma_1=\sigma_2\circ\psi$.  Such	a	pair 	$(\g,\sigma)$	is	said	to be {\sl solvable} (resp.	{\sl simple,...}) if $\g$ is a solvable (resp. simple, ...)
Lie algebra. 
\end{dfn}

\begin{prop} \cite{Lo} There is	a bijection between the set of isomorphism classes of
simply connected symmetric spaces and the set of isomorphism classes of symmetric couples
\end{prop}

Let	us briefly recall how this correspondence is made. If $(M,s )$ is a
symmetric space the group $G=G(M)$ generated by	products of	an even number of symmetries
is a transitive Lie transformation group of $M$, called	the {\sl transvection
group}	of	the symmetric space. One associates to	$(M,s )$ a symmetric couple $(\g,\sigma)$ as follows. The Lie algebra	$\g$ is	the Lie
algebra	 of	the transvection group	$G$.	 Let us choose a	base point $o$ in
$M$ and let $\tilde\sigma$	be	the involutive automorphism of	$G$ obtained
by ``conjugation"	by $s_o$. Then $\sigma$  is the differential of $\tilde\sigma$
at the neutral element	 $e$	of 		$G$. 

\noindent On such a symmetric space, one has a preferred  affine connection.
\begin{prop}\cite{Bith,BCG}
Let $X,Y$ be smooth vector fields on $M$ and let $x$ be any point in $M$. Then the formula:
\begin{equation*}
\left(\nabla_XY\right)_x\;:=\;\frac{1}{2}[X,Y+s_{x\star}Y]_x
\end{equation*}
defines a affine connection $\nabla$ on $M$. The latter has no torsion and its curvature tensor is parallel.
It is called the canonical {\sl Loos connection} on $(M,s)$.
Moreover, the symmetries $\{s_x\}$ coincide with the geoedesic symmetries  w.r.t. the Loos connection.
\end{prop}
The subalgebra $\k$ of $\g$ is known to be isomorphic to
the holonomy algebra for the Loos connection $\nabla$ on the symmetric space
$M=G/_{\textstyle{K}}$---where $K$ is the isotropy at $o$ (its Lie algebra is
isomorphic to $\k$).

\vspace{2mm}

\noindent Now, we recall certain results obtained in \cite{Bi} concerning quantization of symmetric spaces\footnote{The study of the deformed multiplications under oscillatory integral form on symplectic symmetric spaces has been motivated by Weinstein
in \cite{W}.}.

\begin{dfn}
Let $M$ be a symmetric space. We denote by $\mbox{\rm CP}^q(M)$ the space of symmetry-invariant complex-valued smooth functions on $M^q$ ($q\in\N_0$). To every element $c\in\mbox{\rm CP}^q(M)$, we associate
its {\sl coboundary} $\delta c\in\mbox{\rm CP}^{q+1}(M)$:
\begin{equation*}
\delta c(x_0,...,x_q)\;:=\;\sum_{j=0}^q\,(-1)^j\,c(x_0,...,\hat{x_j},...,x_q)\;.
\end{equation*}
One then verifies that $\delta^2=0$.
\end{dfn}
We define equivalently $\delta_{\mbox{\rm op}} :=-\sigma_{12}\circ\delta\circ\sigma_{12}$ where $(\sigma_{12}c)(x_0,x_1,x_2,...,x_q):=
c(x_1,x_0,x_2, x_3,..., x_q)$ and $\sigma_{12}|_{\mbox{\rm CP}^{1}(M)}:=\id$. Note that $\delta_{\mbox{\rm op}}^2=0$ as well.

\noindent Observe that the totally skewsymmetric functions form a  subcomplex of $(\mbox{\rm CP}^\bullet(M),\delta)$.
\begin{dfn}
A totally skewsymmetric function $S\in\mbox{\rm CP}^{3}(M)$ is called {\sl admissible}
if for all $x,y$ and $z$ in $M$, it verifies:
\begin{equation*}
S(x,s_xy,z)\;=\;-S(x,y,z)\;.
\end{equation*}
\end{dfn}
\begin{prop}
Let $M$ be a symmetric space endowed with an invariant volume form. Then  for every admissible function $S\in\mbox{\rm CP}^3(M)$ which is cocyclic (i.e. $\delta S=0$), the following formula:
\begin{equation*}
u\star v (x)\;:=\;\int_{M\times M}e^{S(x,y,z)}\,u(y)\,v(z)\,{\rm d}y\,{\rm d}z\qquad(u,v\in C^\infty_c(M))
\end{equation*}
defines a formally associative product $\star$.
\end{prop}
Analogously to the Lie group case mentioned in the first paragraph, the above statement naturally leads to consider a class of symmetric spaces endowed with symplectic structures\cite{Bith}.
\begin{dfn}
A symplectic  symmetric space is a triple $(M,\omega,s)$ where $(M,s)$ is a symmetric space and where 
$\omega$ is a non-degenerate two-form on $M$ which is invariant under the symmetries.
\end{dfn}
In this situation, the two-form is automatically parallel w.r.t. the Loos connection; in particular $\omega$
turns out to be a symplectic form. At the simply connected level symplectic symmetric spaces are 
encoded by the following algebraic structures particularizing the notion of symmetric couples.
\begin{dfn}
A symplectic triple is a triple $(\g,\sigma,{\bf \Omega})$ where $(\g,\sigma)$ is a symmetric couple
and where ${\bf \Omega}$ is a Chevalley two-cocycle on $\g$ w.r.t. the trivial representation on $\R$
such that its restriction to $\p\times\p$ is non-degenerate (one takes the convention that $i_\k{\bf \Omega}=0$).
\end{dfn}
The non-compact
rank one Hermitean symmetric space $G_s/K_s$ induces, via the Iwasawa decomposition, a left-invariant affine connection on the symplectic Lie group $\S$. Now, there exists other symplectic symmetric spaces that admit
$\S$ as a simply transitive automorphism group.  In particular, there exist other left-invariant  complete symplectic symmetric connections on $\S$ than the above mentioned Kaehler one. Those are simpler to handle in many respects.
They may be viewed as obtained from the Kaehler connection by contracting to zero certain components
(but not all of them!) of its curvature tensor. Hence the terminology `contracted geometry' used in the second 
paragraph. Split solvable Ricci type spaces are of this form.

\section{Star-products on Split solvable Ricci type symplectic symmetric spaces}\label{RICCI}
We start by describing the geometry of the split solvable Ricci type symplectic symmetric spaces
defined in  \cite{CGS} (see also \cite{Bi2} Table (1) Prop. 2.3 for a first occurrence of these spaces in dimension four).

The transvection algebra $\g$ of such a space is a one dimensional split extension of 
two copies of the Heisenberg algebra $\h=V+\R E$ associated to a bilinear symplectic form $\Omega$ on the vector space $V$ :
\begin{equation*}
\g=\a\times(\h\oplus\h)\;;
\end{equation*}
where $\a=\R H$ and $[H,(v,\ell)\oplus(v',\ell')]:=(v,2\ell)\oplus(-v',-2\ell')$ ($(v,\ell),(v',\ell')\in\h$). 
The involution $\sigma$ of $\g$ is defined by
\begin{equation*}
\sigma(aH\oplus X\oplus Y)\;:=\;(-aH)\oplus Y\oplus X\qquad(X,Y\in\h)\;.
\end{equation*}
It yields the ($\pm1$)-eigenspace decomposition:
\begin{equation*}
\g=\k\oplus\p\,,\quad\k:=\Delta\h\subset\h\oplus\h\,\mbox{ and }\p:=\a\oplus\{(X,-X)\}\;,
\end{equation*}
where $\Delta\h$ denotes the diagonal in $\h\oplus\h$.

\noindent The following symplectic structure on $\p$ (within the obvious notations):
\begin{equation*}
{\bf \Omega}\;:=\;{\rm d }a\wedge{\rm d}\ell\,+\,\Omega,
\end{equation*}
is $\k$-invariant. The triple $(\g,\sigma,{\bf \Omega})$ is then  a symplectic triple and we denote
by $(M,\omega, s)$ the corresponding connected simply connected symplectic symmetric space.
Observe that the Lie algebra
\begin{equation*}
\s\;:=\;\a\times\h\;,
\end{equation*}
is supplementary to $\k$ in $\g$. The corresponding analytic subgroup $\S$ of the transvection group
$G=G(M)$ therefore admits an open orbit in $M$.
\begin{prop}\cite{BM,BCSV}
The group $\S$ acts simply transitively on $M$. The map
\begin{equation}\label{DARBOUX}
\s=\a\times V\times\R E\longrightarrow M:(a,v,\ell):=aH\oplus x\oplus \ell E\mapsto\exp(aH)\exp(v)\exp(\ell E) K
\end{equation}
is a global Darboux coordinate system ($omega$ corresponds to ${\bf \Omega}$) within which the geodesic symmetries
read as follows:
\begin{equation*}
\begin{split}
s_{(a,v,\ell)}(a',v',\ell')=(&2a-a',2  \cosh(a-a')v-v',\\
   			&2  \cosh(2 (a -  a'))\ell + \Omega(v,v') \sinh(a - a')-\ell').
\end{split}
\end{equation*}
\end{prop}
\begin{rmk}
{\rm 
In particular, the homogeneous space $M$ is of {\sl group type} meaning that there are subgroups
of automorphisms whose underlying manifolds identify with $M$. Quantizations of $M$, therefore
yield left-invariant star products (formal or not) on these groups---or, {\sl universal deformation formulae} for actions of the latter groups. This aspect is discussed in \cite{BBM, BCSV,Bi07}.
}
\end{rmk}
\begin{prop}
\begin{enumerate}
\item[(i)] The derivation algebra of $(M,\omega,s)$ is isomorphic to $\sp(V,\Omega)\times\g$.
\item[(ii)] $\CU(\s)^\k=\CU(\R.E)$\;.
\item[(iii)] $(\s\wedge\s)^\k\simeq\R\oplus V\;.$
\item[(iv)] $H^2_{\mbox{\tiny Chevalley}}(\s,\R)^\k=0$\;.
\end{enumerate}
\end{prop}
\Pf
(i) We consider the linear injections $\h\to\g:X\mapsto X^\p:=(X,-X)$ and $\h\to\g:X\mapsto X^\k:=(X,X)$.
Let $D\in\Der(\g,\sigma,\Omega)$. 
Within these notations, we observe that:
\begin{eqnarray*}
&&DE^\k=D[X^\p,Y^\p]=(\Omega(DX^\p,Y^\p)+\Omega(X^\p,DY^\p))E^\k=0\;,\\
&& [DH,E^\k]=[DH,E^\k]=-2DE^\p\sim E^\p\;,\\
&& [DE^\p,v^\p]+[E^\p+Dv^\p]=0\;\Rightarrow\;[E^\p,Dv^\p]=0\;\Rightarrow\;DV^\p\subset V^\p\;,\\
&& D[H,v^\p]=Dv^\k=[DH,v^\p]+[H,Dv^\p]=Dv^\k+[DH,v^\p]\;\Rightarrow\;[DH,v^\p]=0\;\Rightarrow DH\sim E^\p\;\\ && \Rightarrow DE^\p=0=DH\;.
\end{eqnarray*}
(ii) Observe first that theCampbell-Baker-Haussdorf formula yields $$\exp(v^\k)\exp(w^\p)=\exp(v^\k+w^\p)\exp(\frac{1}{2}\Omega(v,w)E^\p)\;.$$
Therefore, since $X=\frac{1}{2}(X^\k+X^\p)$ for all $X\in\h$, one has 
\begin{equation*}
\varphi(a,x,z)=\exp(aH)\exp(\frac{1}{2}x^\p)\exp(\frac{1}{2}zE^\p).K\;.
\end{equation*}
Consider $Z\in\k$. Then,
\begin{eqnarray*}
Z^\star_{\varphi}&=&\ddto\exp(-tZ^\k_V)\exp(-tZ^\k_\z)\exp(aH)\exp(\frac{1}{2}x^\p)\exp(\frac{1}{2}zE^\p).K=\\
&=& \ddto\exp(aH)\exp(-t\Ad(\exp(-aH))(Z^\k_V+Z^\k_\z))\exp(\frac{1}{2}x^\p)\exp(\frac{1}{2}zE^\p).K=\\
&=& \ddto\exp(aH)\exp(-t(\cosh(a)Z^\k_V+\sinh(a)Z^\p_V+\cosh(2a)Z^\k_\z+\sinh(2a)Z^\p_\z))\exp(\frac{1}{2}x^\p)\times\\
&\times&\exp(\frac{1}{2}zE^\p).K=\\
&=& \ddto\exp(aH)\exp(
-t(\cosh(a)Z^\k_V+\sinh(a)Z^\p_V))\exp(\frac{1}{2}x^\p)\exp(\frac{1}{2}zE^\p-t\sinh(2a)Z^\p_\z).K=\\
&=& \ddto\exp(aH)\exp(
-t\sinh(a)Z^\p_V)\exp(-t\cosh(a)Z^\k_V)\exp(\frac{1}{2}x^\p)\exp(\frac{1}{2}zE^\p-t\sinh(2a)Z^\p_\z).K=\\
&=& \ddto\exp(aH)\exp(
-t\sinh(a)Z^\p_V)\exp(\frac{1}{2}x^\p)\times\\
&\times&\exp(\frac{1}{2}zE^\p-t\sinh(2a)Z^\p_\z-\frac{1}{2}t\cosh(a)\Omega(Z_V,x)E^\p).K=\\
&=& \ddto\exp(aH)\exp(\frac{1}{2}(
-2t\sinh(a)Z^\p_V+x^\p))\times\\
&\times&\exp(\frac{1}{2}(zE^\p-2t\sinh(2a)Z^\p_\z-t\cosh(a)\Omega(Z_V,x)E^\p)).K=\\
&=&-2\sinh(a)Z_V-\cosh(a)\Omega(Z_V,x)E-2\sinh(2a)Z_\z\;.
\end{eqnarray*}
Moreover, one has the following expressions for the left-invariant vector fields:
\begin{eqnarray*}
\tilde{H}&=&\ddto\exp(aH)\exp(x)\exp(zE)\exp(tH)=\ddto\exp((a+t)H)\exp(e^{-t}x)\exp(ze^{-2t}E)=\\
&=&H-x^j\partial_{x^j}-2zE\,;\\
\tilde{v}&=& \ddto\exp(aH)\exp(x+tv)\exp((z+\frac{1}{2}t\Omega(x,v))E)=\partial_v+\frac{1}{2}\Omega(x,v)E\,;\\
\tilde{E}&=&E\,.
\end{eqnarray*}
In particular, 
$$
Z^\star=-2\left(e^aZ_V-\cosh(a)\tilde{Z_V}+\sinh(2a)\tilde{Z_\z}\right)\;.
$$
Now, let $P=\sum_{J,i}\alpha^{J,i}A_JH^i$ be an element of $\CU(\s)$ where $A_J\in\CU(\h)$
and $\alpha^{J,i}\in\R$. Then, $[E^\star,\tilde{P}]=\sum\alpha^{J,i}\tilde{A_J}[E^\star,\tilde{H}^i]$.
Since $[E^\star,\tilde{H}]=-2\cosh(2a)\tilde{E}-2\sinh(2a)\tilde{E}=-2e^{2a}\tilde{E}$, one concludes that
$[E^\star,\tilde{P}]$ is left-invariant only if $\alpha^{J,i}=0$ for $i>0$. Moreover, for all $X\in V$,
one has $[\partial_v,\tilde{X}]=\frac{1}{2}[\tilde{v},\tilde{X}]$.
 Therefore, for $P\in\CU(\h)$, one gets
$[Z^\star,\tilde{P}]\sim(\frac{1}{2}e^a-\cosh(a))[\tilde{Z_V},\tilde{P}]$ which is left-invariant only if $P$ is central
in $\CU(\h)$.

\vspace{3mm}

\noindent (iii) The Lie bracket yields a map $(\s\wedge\s)^\k\to\s^\k$. Considering 
$w:=h^jH\wedge v_j+\lambda  H\wedge E+\beta^{kl}v_k\wedge v_l+e^i E\wedge v_i$ in $(\s\wedge\s)^\k$,
one then observes that $h^j=0$. The latter implies $[w,w]=0$ as well as Poisson compatibility of the different terms
involved. Therefore, writing $w=\lambda({}^\sharp{\bf \Omega}+B+E\wedge v)$ with $v\in V$ and $B\in V\wedge V$ yields
$0=[Z^\star,B+E\wedge v]\sim f(a)[Z_V,B]$ where $f(a)$ is some non zero function of $a$. A small computation
yields $B=0$.

\vspace{3mm}

\noindent (iv)  By dualizing w.r.t. ${\bf \Omega}$, 
the preceding paragraph implies that every $\k$-invariant  Chevalley-2-cocycle is of the form $\lambda \delta E^\star+H^\star\wedge{}^\flat v$ for some $v\in V$. One then obtains the assertion by observing
that $\delta{}^\flat v=H^\star\wedge {}^\flat v$.
\EPf
\begin{cor}\label{SP}
Up to a change of formal parameter, every two symmetry-invariant formal star products on a 
split Ricci-type solvable symplectic symmetric space are equivalent under the action of an invertible element of 
$\CU(\R.\tilde{E})[[\theta]]$. 
\end{cor}
\section{Strict Quantizations}
\subsection{Recalling a first example}
We first recall a construction of strict quantization on the group $\S$ \cite{BM}.
\noindent Let us denote by
\begin{equation*}
\CF_\z u(a,v,\xi)\;:=:\;\hat{u}(a,v,\xi):=\int e^{-i\xi z}u(a,v,\ell)\,{\rm d}\ell
\end{equation*}
the partial Fourier transform in the $z$-variable. Denoting by $\CS$ (resp. $\tilde{\CS}$)
the space of Schwartz functions in the variables $(a,x,z)$ (resp. $(a,v,\xi)$), one
observes that the above partial Fourier transform establishes an isomorphism $\CF_\z:\CS\to\tilde{\CS}$.

\noindent The following diffeomorphisms will play an important role.
\begin{dfn} Denote by $\tilde{\S}:=\{(a,v,\xi)\}$ the space of variables of $\CF_\z$-transformes functions.
The {\bf twisting map} is the smooth one-parameter family of diffeomorphisms defined as
\begin{equation*}
\phi_\theta:\tilde{\S}\to\tilde{\S}:(a,v,\xi)\mapsto(a,\frac{1}{\cosh(\theta\xi)}v,\frac{1}{2\theta}\sinh(2\theta\xi))\;.
\end{equation*}
The partial map $\xi\mapsto\frac{1}{2\theta}\sinh(2\theta\xi)$ will be denoted by $\varphi_\theta\in\mbox{\rm Diff}(\R)$.
\end{dfn}
\begin{lem}
For all $\theta\in\R$, one has 
\begin{enumerate}
\item[(i)] $\phi^\star_\theta\tilde{\CS}\subset\tilde{\CS}$;
\item[(ii)] $(\phi^{-1}_\theta)^\star\tilde{\CS}\subset\tilde{\CS}'$.
\end{enumerate}
\end{lem}
Note that the map
\begin{equation}\label{TAUS}
T^{-1}_{\theta,1}\;:=\;\CF_\z^{-1}\circ\phi^\star_\theta\circ\CF_\z\,:\,\CS\to\CS
\end{equation}
admits a formal asymtotic expansion
\begin{eqnarray*}
{\tilde{T}}^{-1}_{\theta,1}:C^\infty(\S)[[\theta]]\to C^\infty(\S)[[\theta]]\\
{\tilde{T}}^{-1}_{\theta,1}=I+\sum_k\theta^kU_k
\end{eqnarray*}
as a formal power series with coefficients $\{U_k\}$ in the differential operators on $\S$.
Denoting by ${\tilde{T}}_{\theta,1}:C^\infty(\S)[[\theta]]\to C^\infty(\S)[[\theta]]$ its formal inverse,
 denoting by $\tilde{\star}^0_\theta$ the formal Moyal star product on $(\s,{\bf \Omega})$ and identifying
$(\S,\omega)$ with $(\s,{\bf \Omega})$ via the Darboux map (\ref{DARBOUX}), one gets a
formal star product $\tilde{\star}_{\theta, 1}$ on $C^\infty(\S)[[\theta]]$ defined by the following formula:
\begin{equation*}
u\,\tilde{\star}_{\theta, 1}\,v\;:=\;{\tilde{T}}_{\theta,1}\left({\tilde{T}}^{-1}_{\theta,1}(u)\,\tilde{\star}^0_\theta\,{\tilde{T}}^{-1}_{\theta,1}(v)\right)\;.
\end{equation*}
The latter star product turns out to be symmetry invariant and coincides with the formal asymtotic expansion
of the non-formal product $\star_{\theta,1}$ defined below \cite{BBM}.
\begin{thm}\label{BM}\cite{BM}

\noindent (i) Set 
\begin{equation*}
T_{\theta, 1}\;:=\;\CF_\z^{-1}\circ(\phi^{-1})^\star_\theta\circ\CF_\z
\end{equation*} and consider the following tempered distribution space on $\S$
\begin{equation*}
\CE_{\theta,1}\;:=\;T_{\theta, 1}\,\big(\CS\big)\;.
\end{equation*}
Then, for all $\theta\in\R$, one has the inclusion:
\begin{equation*}
\CS\subset\CE_{\theta,1}
\end{equation*}
and   the linear isomorphism:
\begin{equation*}
T^{-1}_{\theta,1}\;:=\;\CF_\z^{-1}\circ\phi^\star_\theta\circ\CF_\z:\CE_{\theta,1}\to\CS
\end{equation*}
extends the map (\ref{TAUS}) above.

\noindent (ii) Denote by $\star^0_\theta$ the Weyl product on $(\s,{\bf \Omega})$. Then the formula
\begin{equation*}
a\,{\star}_{\theta, 1}\,b\;:=\;{{T}}_{\theta,1}\left({{T}}^{-1}_{\theta,1}(a)\,{\star}^0_\theta\,{{T}}^{-1}_{\theta,1}(b)\right)
\end{equation*}
defines an associative algebra $(\CE_{\theta,1}, {\star}_{\theta, 1})$.

\noindent (iii) For $u,v$ in $C^\infty_c(\S)$, the product reads
\begin{equation*}
u\,{\star}_{\theta, 1}\,v (x)\;=\;\frac{1}{\theta^{\dim(\S)}}\int_{\S\times\S} \,{\bf A}_1(x,y,z)\,e^{\frac{i}{\theta}S(x,y,z)}\,
u(y)\,v(z)\,{\rm d}y\,{\rm d}z\;,
\end{equation*}
where 
\begin{eqnarray*}
{\bf A}_1(x,y,z)\;:=\;\cosh(2(a_1-a_2))\,[\cosh(a_2- a_0)\cosh(a_0-a_1)\,]^{\dim\S-2}\mbox{ and}\\
S(x,y,z)\;:=\;S_0\big(\cosh(a_1-a_2)v_0, 
\cosh(a_2-a_0)v_1, \cosh(a_0-a_1)v_2\big)- \cyclic_{0,1,2}\sinh(2(a_0-a_1))\ell_2\;,
\end{eqnarray*}
with  $S_0(v_0,v_1,v_2):=\Omega(v_0,v_1)+\Omega(v_1,v_2)+\Omega(v_2,v_0)$
which denotes the Weyl product phase on $(V,\Omega)$.

\noindent (iv) Both phase $S$ and amplitude ${\bf A}_1$ are invariant under the diagonal action of the symmetries $\{s_x\}_{x\in M}$ on $M\times M\times M$. Moreover, the  formula:
\begin{equation*}
(\,a\,|\,b\,)_{\theta,1}\;:=\;\int_\S T^{-1}_{\theta,1}(a)\overline{T^{-1}_{\theta,1}(b)}
\end{equation*}
defines a symmetry invariant pre-Hilbert structure on $\CE_{\theta,1}$ such that:
\begin{eqnarray*}
(\,a\,{\star}_{\theta, 1}\,c\,|\,b\,)_{\theta,1}=(\,a\,|\,b\,{\star}_{\theta, 1}\,\overline{c}\,)_{\theta,1}\;.
\end{eqnarray*}

\noindent (v) The product ${\star}_{\theta, 1}$ extends by continuity to the Hilbert completion $\CH_{\theta,1}$ of $\CE_{\theta,1}$ w.r.t. the above pre-Hilbert structure. One then obtain a field of associative Hilbert algebras
on which the symmetries act by unitary automorphisms.
\end{thm}
\subsection{Schwartz multipliers and tracial products}
Corollary \ref{SP} tells us that at the formal level every symmetry invariant star product on $M$ may be obtained
by transporting the above star product $\tilde{\star}_{\theta, 1}$ via a $G$-equivariant equivalence of the 
form
\begin{equation*}
P=\sum_k\theta^k\tilde{P}_k\;\quad\mbox{ with }\tilde{P}_k\in\C[\tilde{E}]\;.
\end{equation*}
That is, denoting by $\CM_f$ the multiplication operator by the function $f$, by an equivalence of the form:
\begin{equation*}
P\;=\;\CF^{-1}_\z\circ\CM_{\CP}\circ\CF_\z\;,
\end{equation*}
where $\CP$ is some formal function on $\tilde{\S}$ depending only on the variable $\xi$: $\CP\in C^\infty_{(\xi)}(\R)[[\theta]]$. Now, a formal computation yields the following `oscillatory integral' formula
for the transported product $\tilde{\star}_{\theta,\CP}$:
\begin{equation*}
u\tilde{\star}_{\theta,\CP} v(x_0)\;=\;\int_{\S\times\S}\frac{(\varphi^\star_\theta\CP)_{\frac{1}{\theta}}(a_1-a_2)}{(\varphi^\star_\theta\CP)_{\frac{1}{\theta}}(a_0-a_2)(\varphi^\star_\theta\CP)_{\frac{1}{\theta}}(a_1-a_0)}\,
{\bf A}_1(x_0,x_1,x_2)\,e^{\frac{i}{\theta}S(x_0,x_1,x_2)}\,u(x_1)\,v(x_2)\,{\rm d}x_1\,{\rm d}x_2\;;
\end{equation*}
where we set $f_\lambda(t):=f(\lambda t)$. In what follows we will render rigourous the above heuristic formula.

\noindent We denote by $\h$ the $G$-invariant foliation in $M$ integrating the distribution associated to $\h^\p$ in $\p$.
\begin{dfn}
We define the space $\mbox{\rm CP}^q(M)_\h$ of $q$-cochains in $\mbox{\rm CP}^q(M)$ that are 
$\h\times\h\times...\times\h$ leafwise constant.
\end{dfn}
\begin{lem} The Darboux map (\ref{DARBOUX}) yields an isomorphism:
\begin{equation*}
C^\infty(\a)\;\simeq\mbox{\rm CP}^2(M)_\h\;.
\end{equation*}
\end{lem}
\Pf
On the the first hand, by invariance, we get $\mbox{\rm CP}^2(M)\simeq C^\infty(M)$.
On the other hand, the quotient  $M\to M/\h$ is a trivial fibration over $M/\h\simeq\R$. Indeed, under the identification
$M=\S$, the foliation leaves correspond to the lateral classes of the normal subgroup $\exp(\h)$ of $\S$.
\EPf
\begin{dfn}
We denote by ${}\bf \Theta$ the subspace of $C^\infty(\R,\mbox{\rm CP}^2(M)_\h)$ constituted by the elements
$\tau:\theta\mapsto\tau_\theta$ such that 

(i) for all $\theta$, the function $\exp\circ\,\tau_\theta$ belongs to the
space $\CO_M(\a)$ of Schwartz multipliers on $\CS(\a)$.

(ii) $(\phi_\theta^\star\tau_\theta)_{\frac{1}{\theta}}\big|_{\theta=0}\equiv0$\;.
\end{dfn}
Note that every element $\tau\in{\bf \Theta}$ defines a linear injection:
\begin{equation*}
T_{\theta,\tau}\;:=\;\CF_\z^{-1}\circ\CM_{\exp(\tau_\theta)}\circ(\phi^{-1}_\theta)^\star\circ\CF_\z\;:\CS\to\CS'\;.
\end{equation*}
We consider the corresponding distribution space:
\begin{equation*}
\CE_{\theta,\tau}\;:=\;T_{\theta,\tau}\big(\CS\big)\;,
\end{equation*}
as well as the linear isomorphism:
\begin{equation*}
T_{\theta,\tau}^{-1}:\CE_{\theta,\tau}\to\CS\;.
\end{equation*}
\begin{thm} 
\noindent (i) For every $\tau\in{\bf \Theta}$, one has the inclusion
\begin{equation*}
\CS\subset\CE_{\theta,\tau}\;.
\end{equation*}

\noindent (ii) On $C^\infty_c(M)\subset\CE_{\theta,\tau}$, the $T_{\theta,\tau}$-transported Weyl product $\star_{\theta,\tau}$ reads
\begin{equation*}
u\star_{\theta,\tau}v\;=\;\frac{1}{\theta^{\dim M}}\int_{M\times M}{\bf A}_1\,\exp\left[\frac{i}{\theta}S+\delta_{\mbox{\rm op}}\left(\sigma_{12}(\phi^\star_\theta\tau_\theta)_{\frac{1}{\theta}}\right)\right]\,u\otimes v\;.
\end{equation*}

\noindent (iii) A symmetry equivariant field of associative Hilbert algebras $\{\CH_{\theta,\tau}\}$ is defined in the same way as in items (vi)-(v) of Theorem \ref{BM}.

\noindent (iv) One has $\CH_{\theta,\tau}\,=\,L^2(M)$ if and only if $|\exp(\tau)|^2=|\mbox{\rm Jac}_{\phi_\theta^{-1}}|$, in which case, the map $T_{\theta,\tau}$ extends as a unitary automorphism of $L^2(M)$.
At the level of compactly supported functions, the corresponding tracial products read:
\begin{equation*}
u\star^{\mbox{\rm tr}}_{\theta,\psi}v\;=\;\frac{1}{\theta^{\dim M}}\int_{M\times M}{\bf A}_{\mbox{\rm can}}\,\exp\left[\frac{i}{\theta}S+i\delta_{\mbox{\rm op}}\left(\sigma_{12}(\phi^\star_\theta\psi_\theta)_{\frac{1}{\theta}}\right)\right]\,u\otimes v\;;
\end{equation*}
where $\psi_\theta$ is arbitrary real-valued and where
\begin{eqnarray*}
{\bf A}_{\mbox{\rm can}}(x_0,x_1,x_2)&:=&
\sqrt{\cosh(2(a_0-a_1))\cosh(2(a_1-a_2))\cosh(2(a_2-a_0))}\,\times \\ &\times&[\cosh(a_0-a_1)\cosh(a_1- a_2)\cosh(a_2- a_0)\,]^{\frac{\dim\S-2}{2}}\;.
\end{eqnarray*}
\end{thm}
\Pf
The product formula in item (ii) is obtained by a straightforward computation of $\CF_\z^{-1}\circ\CM_{\exp(\tau_\theta)}\circ\CF_\z(\star_{\theta,1})$. The NSC of item (iv) and the formula for ${\bf A}_{\mbox{\rm can}}$ follow from the fact that
\begin{equation*}
\int T^{-1}_{\theta, \tau}(u)\overline{T^{-1}_{\theta, \tau}(v)}=\int \phi^\star_\theta(\exp(\tau_\theta)\hat{u}\,\overline{\exp(\tau_\theta)\,\hat{v}})=\int|\mbox{\rm Jac}_{\phi^{-1}_\theta}|\,
|\exp(\tau_\theta)|^2\,\hat{u}\,\overline{\hat{v}}\;.
\end{equation*}
\EPf
\begin{rmk}
{\rm 
Observe the non-cocyclicity of the logarithm of the amplitude in every case.
}
\end{rmk}
\subsection{The quantization space}
\begin{dfn}
For every $\tau\in
{\bf \Theta}$, we call the corresponding product $\star_{\theta,\tau}$ (on
$\CE_{\theta,\tau}$) a {\sl tempered} deformation quantization of $M$.
\end{dfn}
\begin{thm}
Every symmetry invariant formal star product on $M$ is an asymptotic expansion
of a tempered deformation quantization.
\end{thm}
\Pf
Consider the space $\CO_M(\a)$ endowed with its natural Fr\'echet topology \cite{Je}. 
By application of an adapted form of Borel's lemma \cite{BB}, every
element $\tilde{T}$ in $\CU(\R.E)[[\theta]]$ is the Taylor series of a smooth
function $\tau\in C^\infty(\R,\CO_M(\R))$, that is, of an element of ${\bf \Theta}$.
The asymptotic expansion technique developed in \cite{BBM} then yields the assertion.
\EPf
\begin{rmk}
{\rm 
A formal star product on every Ricci type space (not only symmetric) has been recently described
by Cahen, Gutt and Waldmann in \cite{CGW} in the context of Marsden-Weinstein reduction.
It would be interesting to compare their approach to the one presented here.
}
\end{rmk}


\begin{thebibliography}{BGGS04}
\bibitem[Bi95]{Bith}
P.~Bieliavsky, {\sl Espaces sym{\'e}triques symplectiques}, Thesis,
  Universit{\'e} libre de Bruxelles, 1995.  
\bibitem[Bi98]{Bi2}
P.~Bieliavsky, {\sl Four-dimensional simply connected symplectic symmetric spaces},
  Geom. Dedicata \textbf{69} (1998), no.~3, 291--316. 

\bibitem[Bi02]{Bi}
P.~Bieliavsky, {\sl Strict quantization of solvable symmetric spaces}, J. Symplectic
  Geom. \textbf{1} (2002), no.~2, 269--320. 
\bibitem[Bi07]{Bi07} Pierre Bieliavsky; {\sl Deformation quantization for actions of the affine group};
 arXiv:0709.1110\;.

\bibitem[BM01]{BM}
P.~Bieliavsky and M.~Massar, \emph{Oscillatory integral formulae for
  left-invariant star products on a class of {L}ie groups}, Lett. Math. Phys.
  \textbf{58} (2001), no.~2, 115--128. 

\bibitem[BB03]{BB} Bieliavsky, Pierre; Bonneau, Philippe On the geometry of the characteristic class of a star product on a symplectic manifold.  Rev. Math. Phys.  15  (2003),  no. 2, 199--215.

\bibitem[BBM05]{BBM} Bieliavsky, Pierre; Bonneau, Philippe; Maeda, Yoshiaki Universal deformation formulae, symplectic Lie groups and symmetric spaces.  Pacific J. Math.  230  (2007),  no. 1, 41--57. \bibitem[BCG95]{BCG}
P.~Bieliavsky, M.~Cahen, and S.~Gutt, \emph{Symmetric symplectic manifolds and
  deformation quantization}, Modern group theoretical methods in physics
  (Paris, 1995), Math. Phys. Stud., vol.~18, Kluwer Acad. Publ., Dordrecht,
  1995, pp.~63--73.
\bibitem[BCSV07]{BCSV} Pierre Bieliavsky, Laurent Claessens, Daniel Sternheimer, Yannick Voglaire; {\sl
 Quantized anti de Sitter spaces and non-formal deformation quantizations of symplectic symmetric spaces};
arXiv:0705.4179 (to appear in Contemp. Math.,  AMS)\;.
\bibitem[BC99]{CB}
Bourgeois, F.; Cahen, M. A variational principle for symplectic connections.  J. Geom. Phys.  30  (1999),  no. 3, 233--265.
\bibitem[CGS05]{CGS} Cahen, Michel; Gutt, Simone; Schwachh\"ofer, Lorenz Construction of Ricci-type connections by reduction and induction.  The breadth of symplectic and Poisson geometry,  41--57, Progr. Math., 232, BirkhŠuser Boston, Boston, MA, 2005.
\bibitem[CGW07]{CGW}
Michel Cahen, Simone Gutt, Stefan Waldmann;
 {\sl Symplectic Connections of Ricci Type and Star Products}; arXiv:0710.1477\;.

\bibitem[GZ98]{GZ}
A.~Giaquinto and James~J. Zhang, \emph{Bialgebra actions, twists, and universal
  deformation formulas}, J. Pure Appl. Algebra \textbf{128} (1998), no.~2,
  133--151. 
\bibitem[Je89]{Je} Jenkins, Joe W. A characterization of bi-invariant Schwartz space multipliers on nilpotent Lie groups.  Studia Math.  92  (1989),  no. 2, 101--129.
 \bibitem[Loo69]{Lo}
O.~Loos, \emph{Symmetric spaces. {I}: {G}eneral theory}, W. A. Benjamin, Inc.,
  New York-Amsterdam, 1969. 
\bibitem[Ri93]{Ri}
M.~Rieffel, \emph{Deformation quantization for actions of {${\bf R}\sp d$}},
  Mem. Amer. Math. Soc. \textbf{106} (1993), no.~506
\bibitem[We94]{W} A. Weinstein, {\it Traces and triangles in symmetric symplectic 
spaces}, Symplectic geometry and quantization (Sanda and Yokohama, 
1993), Contemp. Math. {\bf 179} (1994), Amer. Math. Soc., Providence, RI, 
261--270. 

\end{thebibliography}
\end{document}